\documentclass[11pt,twoside,a4paper,notitlepage]{article}

\usepackage{amsmath}

\usepackage{amssymb}
\usepackage{amsfonts}
\usepackage[a4paper,margin=2cm,vmargin={1cm,2cm},includeheadfoot]{geometry}

\usepackage{amsthm}
\usepackage{enumerate}
\usepackage[scaled=0.92]{helvet}
\usepackage{srcltx}
\usepackage{graphicx}

\usepackage[labelfont=bf,margin=1cm,justification=justified,labelsep=period]{caption}

\theoremstyle{plain}
\newtheorem{theorem}{Theorem}[section]
\newtheorem{proposition}[theorem]{Proposition}

\newtheorem{lemma}[theorem]{Lemma}

{\theoremstyle{definition}
\newtheorem{remark}[theorem]{Remark}}

{\theoremstyle{definition}
\newtheorem{remarks}[theorem]{Remarks}}

\newcommand{\R}{{\mathbb {R}}}
\newcommand{\N}{{\mathbb{N}}}

\newcommand{\E}{{\mathbb{E}}}   

\newcommand{\wt}{\widetilde}    
\newcommand{\1}{{\bf 1}}        
\newcommand{\Ord}{{\mathcal{O}}}    
\newcommand{\ord}{{\text{\scriptsize$\mathcal{O}$}}}    

\begin{document}

\begin{center}
{\LARGE\bf Inversion of analytic characteristic functions and infinite convolutions
of exponential and Laplace densities}\footnote{
        Acknowledgements. The authors were partially
 supported by grant MTM2009-08869 from the
Ministerio de Ciencia e Innovacion and FEDER. A. Ferreiro-Castilla was also supported by a PhD
grant of the Centre de Recerca Matem\`{a}tica.\\ [0.1cm]
\hspace*{12pt}E-mail-addresses: aferreiro@mat.uab.cat (A. Ferreiro-Castilla),
        utzet@mat.uab.cat (F. Utzet).}
 \\ [.5cm]

{\large \bf Albert Ferreiro-Castilla}  \\ 

 {\it Centre de Recerca Matem\`{a}tica, Apartat 50, 08193 Bellaterra
(Barcelona) Spain.} \\ [0.5cm]

{\large \bf Frederic Utzet}
\footnote{Corresponding author: tel +34 935813470, fax +34 935812790.}

{\it Departament de Matem\`{a}tiques, Edifici C, Universitat
Aut\`{o}noma de Barcelona,\\   08193 Bellaterra (Barcelona) Spain.}

\end{center}

\noindent\hrule

\bigskip

\noindent{\bf Abstract.}

We prove that certain quotients of entire functions are characteristic
 functions.
Under some conditions, the probability measure corresponding to
 a characteristic function of that type has a density which can be
 expressed
as a generalized Dirichlet series, which in turn is an infinite linear
 combination of exponential or Laplace densities.
These results are applied to several examples.

\bigskip


\bigskip


\noindent\hrule

\section{Introduction}

There are  many cases in the literature where a characteristic function $\varphi(t)$ of a
probability measure can be written as $\varphi(t)=1/g(it)$ for $t\in \R$, where $g(z)$ is an
entire function of the complex variable $z$. Two important examples are the square of a
Kolmogorov law and the L\'{e}vy area. Specifically,
the characteristic function of the square of a Kolmogorov law is given by
$\varphi_{1}(t)=\sqrt{2it}/\sin(\sqrt{2it})$, and thus we are setting
$g_{1}(z)=\sin(\sqrt{2z})/\sqrt{2z}$. This example was studied by Dugu\'e \cite{Dugue66}, who
determined that its distribution function is the Jacobi theta function
\begin{equation*}
F_1(x)=\vartheta_{4}(x)=1+2\sum_{k=1}^{\infty}(-1)^{k}e^{-\pi^2k^{2}x/2}\ ,\qquad x>0\ .
\end{equation*}

\noindent The   L\'{e}vy  area  was introduced by L\'{e}vy \cite{Levy51} (with a slightly different  parametrization)
as the random variable given
by
\begin{equation*}
\label{definicio}
\int_{0}^{1}W_{1}(t)dW_{2}(t)-\int_{0}^{1}W_{2}(t)dW_{1}(t)\ ,
\end{equation*}
where $\{W_1(t),\ t\ge 0\}$ and $\{W_2(t),\ t\ge 0\}$ are two independent standard Brownian
motions. L\'{e}vy  \cite{Levy51} deduced that its characteristic function is
$\varphi_{2}(t)={\rm sech}(t)$. In this case we  take $g_{2}(z)=\cos(z)$. The density was also
computed by
L\'{e}vy  \cite{Levy51} and is
\begin{equation*}
f_2(x)=\frac{1}{2}\,{\rm sech}\left(\frac{\pi}{2}\,x\right)
=\frac{e^{-\pi \vert x\vert /2}}{1+e^{-\pi \vert x\vert}}
=\sum_{k=1}^\infty (-1)^{k+1}e^{-\frac{(2k-1)\pi}{2} |x|}\ ,\qquad x\neq0\ .
\end{equation*}

What these two examples share is that regardless of the behaviour of the entire functions $g_1(z)$
and $g_2(z)$, their   inverses are
 characteristic functions. In fact, it is well known that for an entire function $g(z)$ of order
$\rho< 2$ which has only real zeros and $g(0)=1$, the inverse $\varphi(t)=1/g(it)$ is a
characteristic function (see Lukacs \cite[pp. 88 and 212]{Luk70} for an equivalent formulation).
This result follows from Hadamard factorization Theorem (see, for example, Levin \cite[p.
26]{Levin96}) which states that an entire function  of order $\rho<2$ can be written as
$$g(z)=e^{cz}\prod_{n=1}^\infty\left(1-\frac{z}{a_n}\right)e^{z/a_n}\ ,$$
where $\{a_n,\, n\ge 1\}$ are the zeros of $g(z)$ ($c$ must be real so that the characteristic
function assertion is true). Note that such factorization of $g(z)$ induces a factorization of
$\varphi(t)$, and each factor is of the type
\begin{equation*}
\frac{1}{1-it/a_{n}}\,e^{-it/a_{n}}
\end{equation*}

\noindent which is the characteristic function of a translated positive or negative exponential
law. If $\rho<1$, then we can leave aside the exponential part inside the infinite product for a
canonical representation and $c=0$. In any case $\varphi(t)$ is factorized as a convergent
product of characteristic functions, and hence it is a characteristic function due to L\'{e}vy's
continuity Theorem.

The two mentioned examples differ in terms of  the density function. Notice that the density deduced from
$F_{1}(x)$ can be considered  an infinite linear combination of exponential densities, while
$f_{2}(x)$ seems to be a mix of Laplace densities. The difference between both cases lies in the
entire function $g(z)$. On the one hand, $g_{1}(z)$ has order $1/2$ and all its zeros are simple and
positive; as opposed to the second example, where $g_{2}(z)$ has order $1$ and its zeros are
simple but symmetric with respect to the origin.

The aim of this paper is to give some general results suggested by those  and
other  similar examples. It will be shown that when $g(z)$ and
$h(z)$ are both entire functions of order $\rho,\, \rho'\in (0,2)$ respectively, satisfying
$g(0)=h(0)=1$ and with a further condition over their zeros then $\varphi(t)=h(it)/g(it)$ is a
characteristic function. We prove  that when $\rho,\,
\rho'\in (0,1)$, the zeros of $h(z)$ and $g(z)$ are simple and positive, and some other additional
hypotheses, then the probability measure corresponding to $\varphi(t)$ has a density $f(x)$ that can be written as a sum
of exponential type densities, $$f(x)=-\sum_{n=1}^\infty \frac{h(a_n)}{g'(a_n)}\, e^{-a_n x},\
x>0.$$ This series is a generalized Dirichlet series (see Mandelbrojt \cite{Man69}) which  has
very good properties; in particular, we prove that it is uniformly convergent on every compact
subset of $(0,\infty)$, and the cumulative distribution function also turns out to be  a generalized
Dirichlet series.

When the zeros are simple and symmetric, positive and negative, the existence of a density can
be also proved for $\rho=1$ and $h(z)\equiv 1$. Then the density can be written as a series of
Laplace type densities. This case is important because it covers some elements of the double
Wiener chaos (see Janson \cite{Jan97}) as the double It\^o-Wiener integrals where the kernel has
symmetric zeros with multiplicity 2; in particular, the L\'evy area belongs to this  class.

We study some examples. In addition to  the square of a  Kolmogorov law and the L\'evy area, we
consider the law of the first hitting time of a Bessel process, whose characteristic function is
expressed as a quotient of Bessel functions (Kent \cite{Kent78} and Borodin and Salminen
\cite{BS02}). We also show how this technique can be used to invert some Laplace transforms. In
our last example we study a particular case of the Heston model used in mathematical finance and
we prove that the general theory developed in the first part can be applied to it; such a study was
the starting point of this paper.

To summarize, the paper is twofold. On the one hand, it shows a way to construct a rich family of
characteristic functions. On the other, given a characteristic function that can be
identified as belonging to that family, our results provide a procedure to invert that characteristic
function.

\section{Construction of characteristic functions}\label{cons_ch}

In this section we identify a rich family of characteristic functions
which can be constructed by quotients of entire functions.

From now on we will consider that a \emph{strictly increasing sequence of positive numbers},
$\{a_{n}\, n\ge1\}$, is a sequence satisfying $0<a_{1}<a_{2}<\cdots$ and such that
$\lim_{n\to\infty}a_{n}=\infty$.

\begin{proposition}\label{prop_existence}
Consider two entire functions $g(z)$ and $h(z)$ of order lying in $(0,1)$ such that
$g(0)=h(0)=1$, with simple positive zeros given by the strictly increasing sequences of positive numbers
$\{a_n,\, n\ge 1\}$ and $\{b_n,\, n\ge 1\}$ respectively. Assume that $a_n<b_n$ for
all $n$. Then $\varphi(t)=h(it)/g(it)$ is a characteristic function of a probability measure on
$[0,\infty)$.
\end{proposition}

\noindent{\it Proof.}
By Hadamard's factorization Theorem (see, for example, Levin \cite[p. 26]{Levin96}),
$$
g(z)=\prod_{n=1}^\infty\left(1-\frac{z}{a_n}\right)\qquad
    \text{and}
\qquad h(z)=\prod_{n=1}^\infty\left(1-\frac{z}{b_n}\right)\ .
$$ We conclude from the convergence of both products that
\begin{equation}\label{ch_function_inf}
\varphi(t)=\frac{h(it)}{g(it)}
    =\prod_{n=1}^\infty\frac{1-it/b_n}{1-it/a_n}
    =\prod_{n=1}^\infty\left(\frac{a_n}{b_n}
        +\left(1-\frac{a_n}{b_n}\right)\left(1-\frac{it}{a_n}\right)^{-1}\right)\ .
\end{equation}

\noindent Since $0<a_n/b_n<1$ it follows that each factor of the above product is the
characteristic function of the probability measure
$$\frac{a_n}{b_n}\, \delta_0+\left(1-\frac{a_n}{b_n}\right){\mathcal E}{\rm xp}(a_n)\ ,$$
where $\delta_0$ is a Dirac measure at $0$ and ${\mathcal E}{\rm xp}(a_n)$ is an exponential law
with parameter $a_n$. The result is a consequence of L\'{e}vy's continuity Theorem. \quad
$\blacksquare$

\begin{remarks}\label{prop_existence_h1}
\hspace{0pt}
\begin{enumerate}[\bf1)]
\item The condition $g(0)=h(0)=1$ is merely a way to ease the notation, in fact the same result
would be true if we let $g(0)=h(0)\neq0$.

\item By means of standard manipulations, Proposition \ref{prop_existence} holds true for
$h\equiv 1$. In this case each factor of $\varphi(t)$ is the characteristic function of
the probability law ${\mathcal E}{\rm xp}(a_n)$.
\item  A similar result to Proposition \ref{prop_existence} is proved in Corollary 9.16 by Schilling {\it et al.}
\cite{SchSonVon10}, where it is also deduced that the corresponding probability measure belongs to the
 Bondesson class,
and if $h=1$, then the probability measure is in   the class of convolutions of
exponential densities (for the definition of these classes see Schilling {\it et al.}
\cite[Pages 80 and 87]{SchSonVon10}).
\item It is easy to show that the restriction of simple zeros in Proposition
\ref{prop_existence} can be relaxed if we let $a_{n}$ and $b_{n}$ have the same multiplicity for
all $n$.

\end{enumerate}
\end{remarks}

\begin{proposition}\label{prop_existence_kent}
Consider two even entire functions $g(z)$ and $h(z)$ of order lying in $(0,2)$ such that
$g(0)=h(0)=1$. Let $\{a_n,\, n\ge 1\}$ and $\{b_n,\, n\ge 1\}$ be two strictly increasing
sequences of positive numbers; and let $\{\pm a_n,\, n\ge 1\}$ and $\{\pm b_n,\, n\ge 1\}$ be the
simple zeros of $g(z)$ and $h(z)$ respectively; assume that $a_n<b_n$ for all $n$. Then
$\varphi(t)=h(it)/g(it)$ is a characteristic function of a probability measure on $\R$.
\end{proposition}

\noindent{\it Proof.}
By Hadamard's factorization Theorem
$$
g(z)=\prod_{n=1}^\infty\left(1-\frac{z^2}{a_n^2}\right)\qquad
    \text{and}
\qquad h(z)=\prod_{n=1}^\infty\left(1-\frac{z^2}{b_n^2}\right)
$$ due to the symmetry of the entire functions. The rest of the proof follows the
same argument as the derivation of Proposition \ref{prop_existence}. In this setting each factor
of $\varphi(t)=h(it)/g(it)$ is the characteristic function  of the probability measure
$$\frac{a_n^2}{b_n^2}\, \delta_0+\left(1-\frac{a_n^2}{b_n^2}\right){\mathcal L}{\rm aplace}(a_n)\ ,$$
where ${\mathcal L}{\rm aplace}(a_n)$ is a Laplace law with parameter $a_n$. \quad $\blacksquare$

As pointed out in Remark \ref{prop_existence_h1} there is a straightforward generalization of
Proposition \ref{prop_existence_kent} when $h\equiv 1$.

The statements of Propositions \ref{prop_existence} and \ref{prop_existence_kent}
are related to the main results of the following sections. The reader will notice that it is
possible to state more general results using the same ideas and hence with similar derivations.
We may state Proposition \ref{prop_existence_kent} without requesting  the
functions to be symmetric but with their zeros given by the sequences $\{a_n,\, n\ge 1\}$ and
$\{b_n,\, n\ge 1\}$ where  $a_n$ and $b_n$ have the same sign, and $\vert a_n\vert<\vert b_n\vert$.
 In such a case
\begin{equation*}
\varphi(t)=\frac{h(it)}{g(it)}
    =e^{dit}
    \prod_{n=1}^\infty\left(\frac{a_n}{b_n}
        +\left(1-\frac{a_n}{b_n}\right)\left(1-\frac{it}{a_n}\right)^{-1}\right)
    e^{it\left(\frac{1}{b_{n}}-\frac{1}{a_{n}}\right)}
\end{equation*}

\noindent (assuming $d\in \R$) and each factor is the characteristic function of a translation of the probability measure
\begin{equation*}
\left(\frac{a_n}{b_n}\, \delta_0+
\left(1-\frac{a_n}{b_n}\right){\mathcal E}{\rm xp}(a_n)\right)\ .
\end{equation*}

\noindent For example, we can set $\varphi(t)=\text{Ai}(uit)/\text{Ai}(vit)$, where $0<u<v<\infty$
and $\text{Ai}(z)$ is Airy's function. The Airy function is an entire function of order $3/2$
and its zeros are real and negative (see Katori and Tanemura \cite{KT09}). In this case, although
we can state that $\varphi(t)$ is a characteristic function we will not be able to write
 its distribution function explicitly.

\subsection{Finite convolution of exponential and Laplace densities}\label{finite_conv_exp}

The next two lemmas will prove useful for determining the density of the characteristic functions
in Proposition \ref{prop_existence} and \ref{prop_existence_kent}; we will consider a finite
product approximation of $\varphi(t)$ and determine its density so we can take the limit
afterwards. This section and the  one after will develop the first steps of the procedure.

Let us recall some standard notations. Given two probability measures on $(\R,\mathcal{B}(\R))$,
$P_{1}$ and $P_{2}$, we denote by $P_{1}\star P_{2}$ its convolution
\begin{equation*}
P_{1}\star P_{2}(B):=\int_{\R}P_{1}(B-y)P_{2}(dy)\ ,
\end{equation*}

\noindent where $B\in\mathcal{B}(\R)$ and $B-y=\{x-y,\, x\in B\}$. The characteristic function
of $P_{1}\star P_{2}$ is the product of the characteristic functions of $P_{1}$ and $P_{2}$.
Moreover, if $P_{1}$ and $P_{2}$ are absolutely continuous with density $f_{1}$ and $f_{2}$
respectively, then the density of $P_{1}\star P_{2}$ is given by the convolution of $f_{1}$ and
$f_{2}$
\begin{equation*}
f_1 \star f_2(x):=\int_{-\infty}^{+\infty} f_1(y)f_2(x-y)\, dy.
\end{equation*}


\noindent The convolution of $P_1\star \cdots\star P_n$ (resp. $f_1\star \cdots\star f_n$) is
denoted by $\star_{j=1}^n P_j$ (resp. $\star_{j=1}^n f_j$).

%

\bigskip

Lemma \ref{lemma_exp_pos} is well known in the literature (see, for example, problem 12 of
chapter 1 of Feller \cite{Fel66}).

\begin{lemma}\label{lemma_exp_pos}
Fix $0<\lambda_{1}<\lambda_{2}<\ldots <\lambda_{n}$ and define the couple
$A(n):=\prod_{i=1}^{n}\lambda_{i}$ and
$B(k,n):=\prod_{\begin{smallmatrix}i=1\\i\neq k\end{smallmatrix}}^{n}(\lambda_{k}-\lambda_{i})$.
Then $\star_{j=1}^n {\mathcal E}{\rm xp}(\lambda_j)$ has density given by
\begin{equation*}
f_{n}(x)=(-1)^{n+1}A(n)\sum_{i=1}^{n}\frac{e^{-\lambda_{i}x}}{B(i,n)}\qquad x\geq0\ .
\end{equation*}

\end{lemma}

%

\bigskip

The proof of the following result uses an interesting property given by Bondesson (see
\cite{Bondes92}) in the context of generalized gamma convolutions that, in our setup,  states
that if $Y$ is a non--negative random variable with moment generating function $M_Y(u)$, for $u$
in a neighborhood of $0$, and $T$ is a centered normal random variable with variance $2$,
independent of $Y$, then the random variable $X:=\sqrt{Y}\, T$ has moment generating function
$M_X(u)=M_Y(u^2)$. The proof is completed in the line
$$
M_X(u)=\E\big[e^{u\sqrt{Y}\, T}\big]
    =\E\Big[\E\big[e^{u\sqrt{Y}\,T}/Y\big]\Big]=\E\big[e^{u^2\, Y}\big]=M_Y(u^2)\ .
$$

\noindent In fact, the next result is a consequence of the previous lemma and  will be useful for
deriving the density of the characteristic function of Proposition \ref{prop_existence_kent}.

\begin{lemma}\label{lemma_lap_pos}
Fix $0<\lambda_{1}<\lambda_{2}<\ldots <\lambda_{n}$ and define
$E(k,n)=\prod_{\begin{smallmatrix}i=1\\i\neq k\end{smallmatrix}}^{n}
(\lambda_{k}^{2}-\lambda_{i}^{2})$. Then
$\star_{j=1}^n {\mathcal L}{\rm aplace}(\lambda_j)$ has density given by
\begin{equation*}
f_{n}(x)=\frac{(-1)^{n+1}}{2}A^{2}(n)
    \sum_{i=1}^{n}\frac{e^{-\lambda_{i}|x|}}{\lambda_{i}E(i,n)}\qquad x\in\R\ .
\end{equation*}

\end{lemma}

\noindent{\it Proof.}
Let us consider the characteristic function
\begin{equation*}
\tilde{\varphi}(t)=\prod_{j=1}^{n}\left(1-\frac{it}{\lambda_{j}^{2}}\right)^{-1}\ ,
\end{equation*}

\noindent which corresponds to a random variable, $Y$, that is the sum of $n$ independent
exponential random variables with parameters $\{\lambda_{j}^{2},\, 1\le j\le n\}$. By the
previous lemma, the density of $Y$ is
\begin{equation*}
f_{Y}(y)=(-1)^{n+1}A^2(n)\sum_{i=1}^{n}\frac{e^{-\lambda^{2}_{i}y}}{E(i,n)}\qquad y\ge0\ .
\end{equation*}

Let $T$ be a centered normal random variable with variance $2$, hence the random variable
$X:=\sqrt{Y}T$ has characteristic function
\begin{equation*}
\varphi(t)=\prod_{j=1}^{n}\left(1+\frac{t^{2}}{\lambda_{j}^{2}}\right)^{-1}\ ,
\end{equation*}

\noindent which corresponds to the sum of $n$ independent Laplace random variables with
parameters $\{\lambda_{j},\, 1\le j\le n\}$. In order to find the density of $X$, consider the
pair $(X,T)$ and compute its marginal density by means of a change of variables to the pair
$(Y,T)$. The Jacobian determinant is $\frac{2|X|}{T^2}$, which means that the change of variables
is an almost everywhere diffeomorphism of the plane. Fix $x\in\R$, then
\begin{eqnarray*}
f_{X}(x)
&=&\1_{x\ge0}\int_{0}^{\infty}f_{(X,T)}(x,t)dt+\1_{x<0}\int_{-\infty}^{0}f_{(X,T)}(x,t)dt\\
&=&\1_{x\ge0}\int_{0}^{\infty}f_{(Y,T)}\left(\frac{x^2}{t^2},t\right)\frac{2|x|}{t^2}dt
    +\1_{x<0}\int_{-\infty}^{0}f_{(Y,T)}\left(\frac{x^2}{t^2},t\right)\frac{2|x|}{t^2}dt\\
&=&\int_{0}^{\infty}f_{Y}\left(\frac{x^2}{t^2}\right)f_{T}(t)\frac{2|x|}{t^2}dt\\
&=&\int_{0}^{\infty}
    \sum_{i=1}^n \frac{(-1)^{n+1}A^{2}(n)}{E(i,n)}e^{-\lambda_{i}^2  x^2/t^2}
    e^{-t^2/4}\frac{1}{2\sqrt{\pi}}\frac{2|x|}{t^2}dt \\
&=&\sum_{i=1}^n \frac{(-1)^{n+1}A^{2}(n)|x|}{\sqrt{\pi}E(i,n)}
    \int_{0}^{\infty}\frac{e^{-\lambda^2_i  x^2/t^2}e^{-t^2/4}}{t^2}dt\\
&=&\frac{(-1)^{n+1}}{2}A^{2}(n)
    \sum_{i=1}^{n}\frac{e^{-\lambda_{i}|x|}}{\lambda_{i}E(i,n)}
\end{eqnarray*}

\noindent as required. \quad $\blacksquare$

\subsection{Finite density approximation of the characteristic function}\label{finite_charac}

We derive here the density function for a finite approximation of the characteristic functions of
Proposition \ref{prop_existence} and \ref{prop_existence_kent}.

\begin{lemma}\label{tech_main}
Let $0<a_{1}<a_{2}<\cdots<a_{n}$ and $0<b_{1}<b_{2}<\cdots<b_{n}$ such that $a_{i}<b_{i}$ for
$1\le i\le n$. Write either
\begin{enumerate}[(a)]
    \item \label{tech_1}
    $g_{n}(z)=\prod_{i=1}^{n}\left(1-z/a_{i}\right)\ $ and
    $\ h_{n}(z)=\prod_{i=1}^{n}\left(1-z/b_{i}\right)$, or
    \item \label{tech_2}
    $g_{n}(z)=\prod_{i=1}^{n}\left(1-z^2/a_{i}^2\right)\ $ and
    $\ h_{n}(z)=\prod_{i=1}^{n}\left(1-z^2/b_{i}^2\right)$.
\end{enumerate}

\noindent Let $\varphi_{n}(t)$ be the characteristic function
$\varphi_{n}(t)=h_{n}(it)/g_{n}(it)$. Then the corresponding law of $\varphi_{n}(t)$ is
\begin{equation*}
\left(\prod_{i=1}^{n}\frac{a_{i}}{b_{i}}\right)\delta_{0}+\mu_{n}\ ,
\end{equation*}

\noindent where $\mu_{n}$ is a finite measure on $(0,\infty)$ -- case (\ref{tech_1}) -- or on
$\R\setminus\{0\}$ -- case (\ref{tech_2}) --, with density given by
\begin{equation*}
\text{(\ref{tech_1})}\quad
\frac{d\mu_{n}}{dx}=-\sum_{i=1}^{n}\frac{h_{n}(a_{i})}{g_{n}'(a_{i})}\,
        e^{-a_{j}x}\quad x>0
\qquad\text{or}\qquad
\text{(\ref{tech_2})}\quad
\frac{d\mu_{n}}{dx}=-\sum_{i=1}^{n}\frac{h_{n}(a_{i})}{g_{n}'(a_{i})}
       \, e^{-a_{j}|x|}\quad x\neq0\ .
\end{equation*}
\end{lemma}

\noindent{\it Proof.}
We will prove the result for the case (\ref{tech_1});  case (\ref{tech_2}) is similar.
Since the convolution of two finite measures is absolutely continuous when one of them is, we
abuse of the notation and write
\begin{equation*}
\left(\frac{a_j}{b_j}\, \delta_0
    +\left(1-\frac{a_j}{b_j}\right){\mathcal E}{\rm xp}(a_j)\right)(x)
\end{equation*}

\noindent to denote the probability \emph{density} function of the measure which is a mixture of
a delta measure and an exponential distribution. Let us denote by $I_{n}^{i}$ the set of all
subsets of $\{1,2,\ldots,n\}$ of cardinal $1\leq i\leq n$. Then
\begin{align*}
\star&_{j=1}^n\left(\frac{a_j}{b_j}\, \delta_0
    +\left(1-\frac{a_j}{b_j}\right){\mathcal E}{\rm xp}(a_j)\right)(x)=\\
&=\left(\prod_{i=1}^{n}\frac{a_{i}}{b_{i}}\right)\delta_{0}(x)
        +\sum_{i=1}^{n}\sum_{J\in I_{n}^{n-i}}
        \Bigg[\prod_{k\in J}\frac{a_{k}}{b_{k}}\Bigg]
        \Bigg[\prod_{k\in J^{c}}\left(1-\frac{a_{k}}{b_{k}}\right)\Bigg]
        \Bigg[\sum_{k\in J^{c}}\frac{a_{k}}
        {\prod_{r\in J^{c}\setminus \{k\}}\left(1-\frac{a_{k}}{a_{r}}\right)}
        e^{-a_{k}x}\Bigg] \\
&=\left(\prod_{i=1}^{n}\frac{a_{i}}{b_{i}}\right)\delta_{0}(x)
    +\sum_{k=1}^{n}\frac{e^{-a_{k}x}a_{k}}
    {\prod_{\begin{smallmatrix}r=1\\r\neq k\end{smallmatrix}}^{n}
        \left(1-\frac{a_{k}}{a_{r}}\right)}
        \sum_{i=1}^{n}
        \sum_{\begin{smallmatrix}J\in I_{n}^{n-i}\\k\notin J\end{smallmatrix}}
        \Bigg[\prod_{r\in J}\frac{a_{r}}{b_{r}}\left(1-\frac{a_{k}}{a_{r}}\right)\Bigg]
        \Bigg[\prod_{r\in J^{c}}\left(1-\frac{a_{r}}{b_{r}}\right)\Bigg]\\
&=\left(\prod_{i=1}^{n}\frac{a_{i}}{b_{i}}\right)\delta_{0}(x)
    +\sum_{k=1}^{n}e^{-a_{k}x}
    \frac{a_{k}}{\prod_{\begin{smallmatrix}r=1\\r\neq k\end{smallmatrix}}^{n}
    \left(1-\frac{a_{k}}{a_{r}}\right)}
    \prod_{r=1}^{n}\left(\left(1-\frac{a_{r}}{b_{r}}\right)+
    \frac{a_{r}}{b_{r}}\left(1-\frac{a_{k}}{a_{r}}\right)\right)\\
&=\left(\prod_{i=1}^{n}\frac{a_{i}}{b_{i}}\right)\delta_{0}(x)
    +\sum_{k=1}^{n}e^{-a_{k}x}
    \frac{\prod_{r=1}^{n}\left(1-\frac{a_{k}}{b_{r}}\right)}
    {\frac{1}{a_{k}}\prod_{\begin{smallmatrix}r=1\\r\neq k\end{smallmatrix}}^{n}
    \left(1-\frac{a_{k}}{a_{r}}\right)}\ ,
\end{align*}

\noindent where we have used Lemma \ref{lemma_exp_pos}.\qquad $\blacksquare$

\begin{remark}\label{tech_main_h1}
If $h\equiv1$ in Lemma \ref{tech_main} then the law of $\varphi_{n}(t)$ is just $\mu_{n}$.
\end{remark}

Lemma \ref{tech_main} shows that the distribution has an atom at the origin. This causes a handicap
towards applying the conventional limit theorems to the expressions of the densities therein.
Section \ref{existence_d} will give further assumptions and criteria to overcome this difficulty.
Finally, notice that the key point to this end is  determining the behaviour of
\begin{equation*}
\lim_{n\to\infty}\frac{d\mu_{n}}{dx}\ .
\end{equation*}

\section{Fundamental lemmas}

The following lemmas are essential for the purposes of the paper and they give the technical results
to allow us to use the standard convergence results on the expressions of Lemma \ref{tech_main}.
Before that, let us  remark on a key property of an entire function: the order of an entire function is
greater or equal to the exponent of convergence of its zeros (see Titchmarsh
\cite[p. 251]{Tit52}). That means, let $f$ be an entire function of order $\rho$ and
$\{a_{n},\,n\ge1\}$ its zeros, then
\begin{equation*}
\rho\ge\inf\{\alpha>0\,:\, \sum_{n=1}^{\infty}|a_{n}|^{-\alpha}<\infty\}\ .
\end{equation*}

\noindent In particular, $\sum_{n=1}^{\infty}|a_{n}|^{-\beta}<\infty$ for $\beta>\rho$. The
equality holds if $f$ has a representation as a canonical product, in fact Proposition
\ref{prop_existence} and \ref{prop_existence_kent} give such representation for $h$ and $g$.

\begin{lemma}\label{fundamental}
Let $g(z)$ be an entire function of order $\rho\in (0,1)$ such that $g(0)=1$, with simple positive
zeros given by the strictly increasing sequence of positive numbers $\{a_n,\, n\ge 1\}$.
 Let $h(z)$ be another
function that either
\begin{enumerate}[(a)]
    \item $h(z) \equiv 1$ or

    \item $h(z)$ is an entire function of order $\rho'\in (0,1)$ such that $h(a_{n})\neq0$
    for all $n$.
\end{enumerate}

\noindent Then, for every $x>0$, the series $\sum_{n\ge1}\frac{h(a_n)}{g'(a_n)}\, e^{-a_n x}$
is convergent.
\end{lemma}

\begin{figure}[htb]
\begin{minipage}[t]{0.45\linewidth} 
\centering
\includegraphics[width=0.9\linewidth, height=0.9\linewidth]{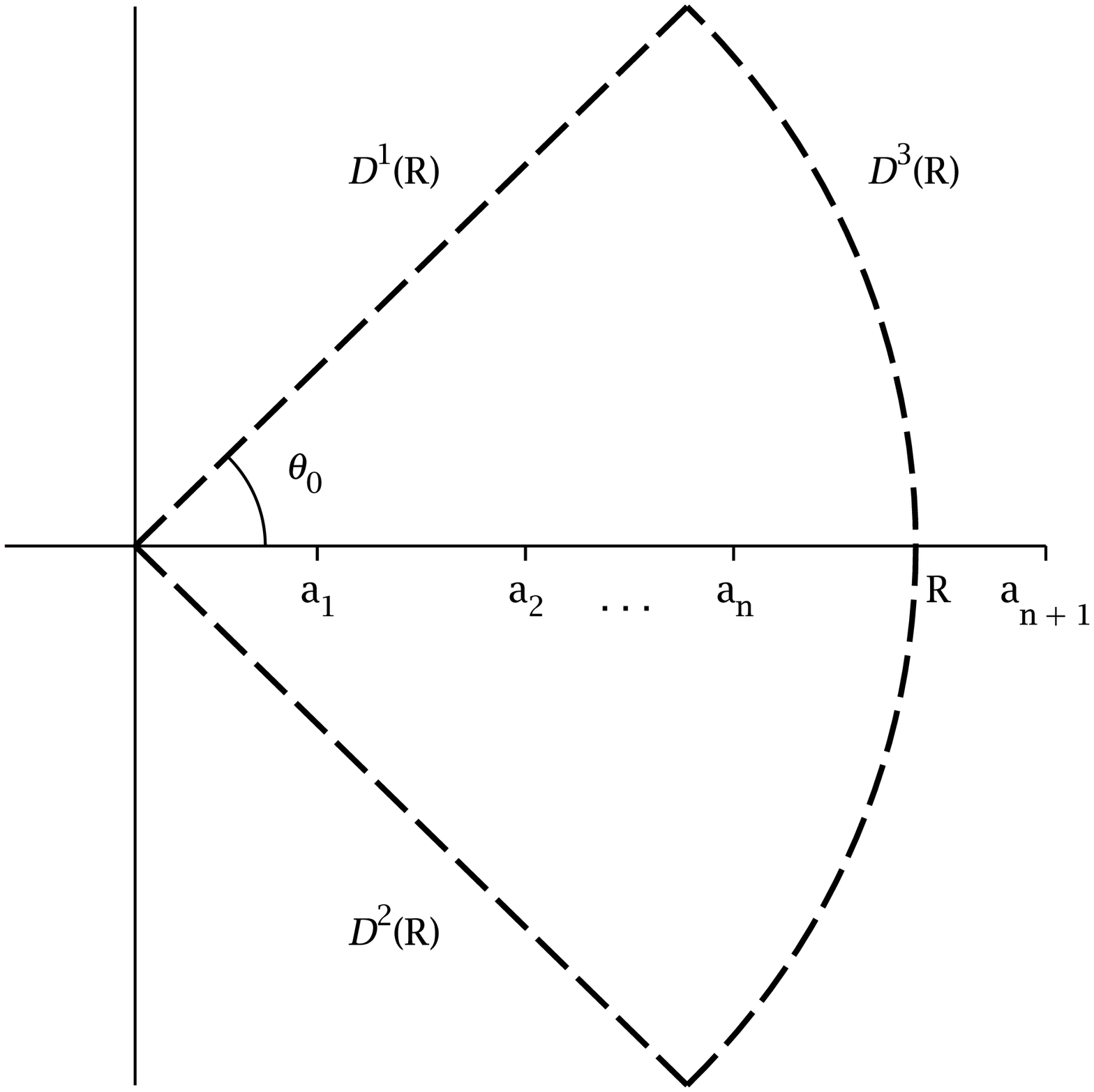}
\caption{Dashed line is the contour $D(R)$ of integration.}
\label{contorn}
\end{minipage}
\hspace{0.05\linewidth} 
\begin{minipage}[t]{0.45\linewidth}
\centering
\includegraphics[width=0.9\linewidth, height=0.9\linewidth]{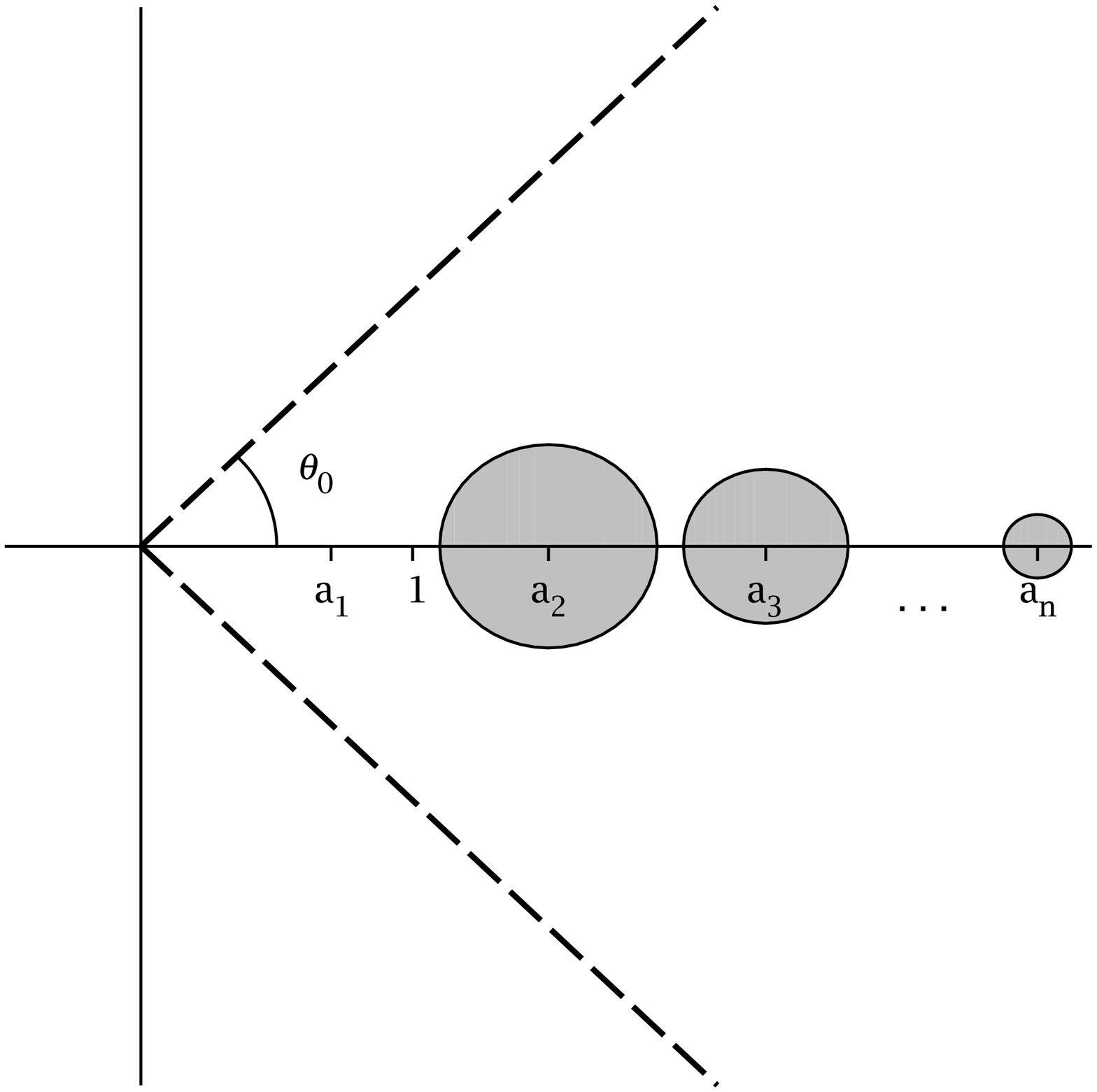}
\caption{The function $\vert g(z) \vert$ is bounded below in the region excluded from the gray
circles.}
\label{excluded}
\end{minipage}
\end{figure}

\noindent{\it Proof.}
We will prove the lemma for  case (b);  case (a) is similar and indeed easier. Consider the
closed  contour $D(R)$ in Figure \ref{contorn}, where $R>0$, $\theta_0\in(0,\pi/2)$ and
$R\ne a_j$ for all $j$.
Fix $x>0$. The function $\frac{h(z)}{g(z)}\, e^{-xz}$ is analytic in the region
bounded by $D(R)$, except at the zeros  $a_j$ that lie in the interior of $D(R)$. The Residue
Theorem allows us to write
$$\frac{1}{2\pi i}\oint_{D(R)} \dfrac{h(z)}{g(z)}\, e^{-xz}\, dz
=\sum_{a_j<R} \text{Res}(a_j)=\sum_{a_j<R}\frac{h(a_j)}{g'(a_j)}\, e^{-a_j x}\ .$$
Our objective is to prove that there is an increasing sequence,
$\{R_n,\,  n\ge 1\}$, towards infinity such that the limit
\begin{equation}\label{finite_lim}
\lim_n \oint_{D(R_{n})} \dfrac{h(z)}{g(z)}\, e^{-xz}\, dz
\end{equation}

\noindent exists and is finite. The key point of the proof is that an entire function of order
$\rho<1$ can be bounded below in the region excluded from  small circles around its zeros. Since the zeros of $g(z)$ are
located on the positive real axis, convenient bounds of $1/|g(z)|$ can be obtained.
Specifically, the proof is based on:
\begin{enumerate}[\bf 1.]
\item \label{fundamental_o_gran}For every $\varepsilon>0,$
    $g(z)=\Ord(e^{\epsilon \vert z \vert})$
    and $h(z)=\Ord(e^{\epsilon \vert z \vert})$. This is due to
    $\sum_n |a_n|^{-1}<\infty$, and
    hence we can apply Titchmarsh's  issue 15  \cite[p. 286]{Tit52}.

\item \label{fundamental_arch}As a consequence of the previous statement, for every
    $\varepsilon>0$ there is an increasing sequence $\{R_n, n\ge 1\}$ such that $\lim_n
    R_n=\infty$ and $\vert g(R_ne^{i\theta}) \vert > \exp\{-\varepsilon R_n\}$ uniformly on
    $\theta\in [0,2\pi]$. See Titchmarsh \cite[p. 276, it. 8.75]{Tit52}.

\item \label{fundamental_rho_p}For every $\varepsilon>0$, there is $\theta_0 \in (0,\pi/2)$
    and $r_{0}>0$ such that
    $$\vert g(re^{\pm i\theta_{0}}) \vert > \exp\{-\varepsilon r\}\ ,
        \ \text{for}\  r\ge r_0.$$
    This is deduced from  Titchmarsh \cite[p. 273, it. 8.71]{Tit52}. Take straight lines
    through the origin with angles $\theta_0$ and $-\theta_0$ that do not intersect with any
    of the circles with centers $a_j$ and radius $1/a_j$  (for $a_j>1$), see Figure
    \ref{excluded}. Titchmarsh proves that for all $\varepsilon'>0$ exists $r_{0}'$ such that
    if $r\ge r_0'$ we have $\vert g(z) \vert > \exp\{- r^{\rho+\varepsilon'}\}$ in the region
    excluded from these discs. Take $\varepsilon'$ satisfying $\rho+\varepsilon'<1$ and thus
    for sufficiently large  $r$ we have $r^{\rho+\varepsilon'}<\varepsilon r$. Choose $r_0$ a
    slightly  larger than $r_{0}'$ in order to obtain the claimed inequality.

\end{enumerate}

Now we are ready to complete the proof. Fix $\varepsilon>0$ such that
$2\varepsilon <x\cos\theta_0$, where $\theta_0$ is the angle depending on $\varepsilon$ such that
property {\bf \ref{fundamental_rho_p}} is fulfilled. Denote by $D^{1}(R_{n})$ and $D^2(R_{n})$
the lower and  upper  straight segment of $D(R_n)$ while $D^3(R_n)$ stands for the  arch. All
paths are considered with the corresponding orientation. Then we divide the integral into three
parts :
$$\oint_{D(R_n)}\frac{h(z)}{g(z)}\,e^{-xz}\,dz=\int_{D^1(R_n)}+\int_{D^2(R_n)}+\int_{D^3(R_n)}.$$

We first consider the integral over $D^3(R_n)$. Due to points {\bf \ref{fundamental_o_gran}}
and {\bf \ref{fundamental_arch}}, we can bound the module of this integral in the following way:
\begin{align*}
\left\vert \int_{D^3(R_n)} \frac{h(z)}{g(z)}\, e^{-xz}\, dz \right\vert & \le R_n
    \int _{-\theta_0}^{\theta_0}\left\vert \frac{h(R_ne^{i\theta})}{g(R_ne^{i\theta})}
    \right\vert\, e^{-xR_n \cos \theta}\, d\theta \\
&\le R_n e^{-xR_{n} \cos \theta_0}  \int _{-\theta_0}^{\theta_0}\left\vert
    \frac{h(R_ne^{i\theta})}{g(R_ne^{i\theta})}  \right\vert \, d\theta
    \le K \theta_0 R_n e^{ -R_n( x \cos \theta_0-2\varepsilon)},
\end{align*}

\noindent for $R_{n}$ large enough. This goes to zero as $n\to \infty$.

The integral over $D^1(R_n)$ can be parametrized as
\begin{equation*}
\int_{D^1(R_n)} \frac{h(z)}{g(z)}\, e^{-xz}\, dz=e^{-i\theta_0}\int_0^{R_n}
\frac{h(re^{-i\theta_0})}{g(re^{-i\theta_0})}\, e^{-xre^{-i\theta_0}} dr\ .
\end{equation*}

\noindent We claim that
\begin{equation}\label{int-inf}
\int_0^{\infty} \left\vert \frac{h(re^{-i\theta_0})}{g(re^{-i\theta_0})}
\, e^{-xre^{-i\theta_0}}\right\vert\, dr<\infty.
\end{equation}

\noindent According to observations {\bf\ref{fundamental_o_gran}} and
{\bf\ref{fundamental_rho_p}}, there is $r_0>0$ depending on $\varepsilon$ such that
$$\int_{r_0}^\infty \left\vert \frac{h(re^{-i\theta_0})}{g(re^{-i\theta_0})}
\, e^{-xre^{-i\theta_0}}\right\vert\, dr
=   \int_{r_0}^\infty \left\vert
    \frac{h(re^{-i\theta_0})}{g(re^{-i\theta_0})} \right\vert\, e^{-xr\cos \theta_0}\,dr
\le K \int_{r_0}^\infty e^{-r(x \cos \theta_0-2\varepsilon)}\,dr<\infty.$$
Since $g(re^{-i\theta _0})\ne 0$ for $r\ge 0$, it turns out that the function
$\left\vert \dfrac{h(re^{-i\theta_0})}{g(re^{-i\theta_0})} \right\vert\, e^{-xr\cos \theta_0}$ is
continuous on $[0,r_0]$ and hence
$$\int_0^{r_0} \left\vert \frac{h(re^{-i\theta_0})}{g(re^{-i\theta_0})}
\, e^{-xre^{-i\theta_0}}\right\vert\, dr<\infty.$$
Adding up the two upper bounds we obtain (\ref{int-inf}). For $D^2(R_n)$ the computations
are equivalent and the limit (\ref{finite_lim}) exists and is finite.\quad $\blacksquare$

\begin{lemma}\label{uniform_absolut}
Under assumptions of Lemma \ref{fundamental}, the series
\begin{equation}\label{general_D}
\sum_{n\geq1}\frac{h(a_n)}{g'(a_n)}\, e^{-a_n x}
\end{equation}

\noindent converges absolutely for $x>0$, and the convergence is uniform on any compact
subset of $(0,\infty)$.
\end{lemma}

\noindent{\it Proof.}
The above expression is a generalized Dirichlet series. As a consequence, if (\ref{general_D}) is
convergent for some $x_0$, then it is so for all $x>x_0$ and the convergence is uniform on every
compact subset of the half-line (see Mandelbrojt \cite[p. 9]{Man69}). Therefore the second part
of the lemma follows from Lemma \ref{fundamental}. Each general Dirichlet series has  associated
an abscissa $\sigma_{c}$ of convergence and an abscissa $\sigma_{a}$ of absolute
convergence. In general these two values are not equal but their distance is bounded (see
Mandelbrojt \cite[p. 11]{Man69}) by
\begin{equation*}
0\le\sigma_{a}-\sigma_{c}\le\limsup_{n\to\infty}\frac{\ln n}{a_{n}}\ .
\end{equation*}

\noindent From the previous observation we know that $\sigma_{c}\le0$; we claim that the above
limit is zero to derive the first part of the result. The claim is proved through Jensen's
formula which gives the relationship $n(r)=\Ord(r^{\rho+\varepsilon})$ for the entire
function $g(z)$, where $\varepsilon>0$ and $n(r)$ stands for the number of zeros with norm less
than or equal to $r$ (see Titchamrsh \cite[p. 249]{Tit52}). Choose $\varepsilon$ such that
$\rho+\varepsilon<1$, then $n\le K a_{n}^{\rho+\varepsilon}$ for some positive constant $K$.
Finally
\begin{equation*}
\ln n\le n-1\le K a_{n}^{\rho+\varepsilon}-1\ ,
\end{equation*}
and hence $\sigma_{c}=\sigma_{a}$.\quad $\blacksquare$

\begin{remark}\label{fundamental_laplace}
Straightforward manipulations lead to generalizing Lemmas \ref{fundamental} and
\ref{uniform_absolut} for the case where the simple zeros of $g(z)$ are
$\{\pm a_{n},\, n\ge1\}$, but still of an  order lying in $(0,1)$.
\end{remark}

The next result somehow extends  case (a) of Lemma \ref{fundamental} but with a penalization of
extra hypotheses. As  will be shown in the examples, this particular case is also of great
interest.

\begin{lemma}\label{fundamental_lap}
Consider an even entire function $g(z)$ of order $1$ such that $g(0)=1$ and
$g(z)=\Ord(e^{A|z|})$ for some positive constant $A$. Let $\{a_n,\, n\ge 1\}$ be a strictly
increasing sequence of positive numbers; let $\{\pm a_n,\, n\ge 1\}$ be the simple zeros of
$g(z)$; assume that  the following limit exists:
\begin{equation}\label{hypo_lap}
\lim_{n\to\infty}\frac{n}{a_{n}}=\delta>0\ .
\end{equation}

\noindent Then, for every $x>0$, the series $\sum_{n\ge1}\frac{1}{g'(a_n)}\, e^{-a_n x}$ is
absolutely convergent and the convergence is uniform on any compact subset of $(0,\infty)$.
\end{lemma}

\noindent{\it Proof.}
It is clear from (\ref{hypo_lap}) that $\sigma_{a}=\sigma_{c}$, thus if we prove the plain
convergence of the series in $(0,\infty)$ the result will follow. The proof is very similar to
that of Lemma \ref{fundamental} but the fact that the order of $g(z)$ is $1$ demands some
modifications. We will consider as before the closed  contour $D(R)$ in Figure \ref{contorn},
where $R>0$, $\theta_0\in(0,\pi/2)$ and $R\ne a_j$ for all $j$. Our objective is to prove that
there is an increasing sequence, $\{R_n,\,  n\ge 1\}$, towards infinity such that
\begin{equation}\label{limit_lap}
\lim_n \oint_{D(R_{n})} \dfrac{1}{g(z)}\, e^{-xz}
\end{equation}

\noindent exists and is finite for a fixed $x>0$. Split the above integral as in
Lemma \ref{fundamental} into $D(R_{n})=D^{1}(R_{n})\cup D^{2}(R_{n})\cup D^{3}(R_{n})$.

The function $g(z)$ admits an expression as a canonical product given by
$$g(z)=\prod_{n=1}^\infty\left(1-\frac{z^2}{a^2_n}\right)\ .$$ Define the function
$$\wt g(z):=\prod_{n=1}^\infty\left(1-\frac{z}{a^2_n}\right)$$
of order $1/2$ and notice $\wt g(z^2)=g(z)$. Fortunately, the function $g(z)$ inherits the good
properties of $\wt g(z)$. To start with, since $\wt g(z)$ has order $1/2$, for every
$\varepsilon>0$ there is a sequence $\wt R_n\nearrow\infty$ such that for all
$\theta\in [0,2\pi]$ $$\vert \wt g(\wt R_ne^{i\theta})\vert >M_{\wt g}(\wt R_n)^{-\epsilon},$$
where $M_{\wt g}(r)=\max_{\{\vert z\vert=r\}}\vert \wt g(z)\vert$ (see Titchmarsh
\cite[p. 275, it. 8.74]{Tit52}). Since $g(z)$ is even, we deduce that there is a sequence
$R_n:=\wt R_{n}^{1/2}\nearrow \infty$ such that for all $\theta\in [0,2\pi]$
$$ \vert g(R_ne^{i\theta})\vert >M_{g}(R_n)^{-\epsilon}\ .$$
Given that $g(z)=\Ord(e^{A\,\vert z\vert})$, for $R_n$ large enough,
$$\left\vert \frac{1}{g(R_ne^{i\theta})}\right\vert
    \le C e^{\varepsilon A R_n} \quad \text{uniformly on } \theta\in [0,2\pi].$$
Now consider $\varepsilon$ such that $A\varepsilon<x\cos\theta_{0}$, and the bound
\begin{align*}
\left\vert \int_{D^3(R_n)} \frac{1}{g(z)}\, e^{-xz}\, dz \right\vert & \le R_n
    \int _{-\theta_0}^{\theta_0}\left\vert \frac{1}{g(R_ne^{i\theta})}
    \right\vert\, e^{-xR_n \cos \theta}\, d\theta \\
&\le R_n e^{-xR_{n} \cos \theta_0}  \int _{-\theta_0}^{\theta_0}\left\vert
    \frac{1}{g(R_ne^{i\theta})}  \right\vert \, d\theta
    \le K \theta_0 R_n e^{ -R_n( x \cos \theta_0-A\varepsilon)},
\end{align*}

\noindent for $R_{n}$ large enough. This goes to zero as $n\to \infty$.

In order to bound the integral over $D^{1}(R_{n})$ consider the following Remark in
Levin \cite[p. 82, eqn. $(2')$]{Levin96}:
\begin{equation*}
\ln|\wt g(\wt re^{i\wt \theta})|= \wt r^{1/2}\pi\delta\sin(\wt \theta/2)
+\frac{\ord(r^{1/2})}{\sin(\wt \theta/2)}
\end{equation*}
for $\wt \theta\in(0,2\pi)$ and as $\wt r\to\infty$; that asymptotic result is based on
the fact that $\widetilde g$ is a canonical product of order 1/2.  Notice that
$\wt g(\wt re^{i\wt \theta})=g(\wt r^{1/2}e^{i\wt \theta/2})
    :=g(re^{i\theta})=g(re^{i(\theta+\pi)})$,
hence the above asymptotic equation is translated into $g(z)$ as
\begin{equation*}
\ln|g(re^{i\theta})|= r\pi\delta|\sin(\theta)|+\frac{\ord(r)}{|\sin(\theta)|}
\end{equation*}
for $\theta\in(0,2\pi)\setminus\{\pi\}$. Therefore
\begin{equation*}
\left|\int_{D^1(R_n)} \frac{1}{g(z)}\, e^{-xz}\, dz\right|
\le
\int_0^{\infty} \left\vert \frac{1}{g(re^{-i\theta_0})}\right\vert
\, e^{-xr\cos\theta_0}\, dr
\sim
\int_0^{\infty}
e^{-r\left(
    x\cos\theta_0
    +\pi\delta|\sin(\theta_0)|
    +\frac{\ord(r)}{r|\sin(\theta_0)|}
    \right)}\, dr\ ,
\end{equation*}
where $\sim$ means that both converge or diverge together. Clearly, the last integral is finite. The
same derivation is valid for the integral over $D^2(R_n)$ and (\ref{limit_lap}) holds.
\quad $\blacksquare$

\section{Construction of density functions}\label{main_S}

This section will derive the density function associated to the characteristic functions of
Proposition \ref{prop_existence} and Proposition \ref{prop_existence_kent} under the
restrictions of Lemma \ref{fundamental_lap}.

Lemma \ref{tech_main} shows the density function for a finite product approximation of the
characteristic function. From L\'{e}vy's continuity theorem, this means a convergence in
distribution of the associated probability measures. We first show a pointwise convergence of the distribution
function to finally obtain the convergence of densities. In order to avoid  Dirac's delta
measure at $0$, instead of using the distribution function $F(x)$ we work with
$\overline{F}(x):=1-F(x)$.

\begin{theorem}\label{main}
Under assumptions of Proposition \ref{prop_existence}, the probability measure on $[0,\infty)$
corresponding to the characteristic function $\varphi(t)=h(it)/g(it)$ is absolutely continuous on
$(0,\infty)$ with (perhaps defective) density given by
$$
f(x)=-\sum_{n=1}^\infty \frac{h(a_n)}{g'(a_n)}\, e^{-a_n x},\ \quad x>0\ .
$$
\end{theorem}

\noindent{\it Proof.}
Recall $h_{n}(z)$, $g_{n}(z)$ and $\varphi_{n}(t)$ from Lemma \ref{tech_main} case
(\ref{tech_1}). Denote by $F_{n}(x)$ the distribution function corresponding to the
characteristic function $\varphi_{n}(t)$. It is clear that $\varphi_{n}(t)\to\varphi(t)$
pointwise, where $\varphi(t)$ is defined in (\ref{ch_function_inf}). Since $a_{n}\neq0$ for all
$n$, it follows that  $\varphi(t)$ is continuous at $0$. Therefore there exists a distribution
function $F(x)$, such that
\begin{equation*}
\lim_{n\to\infty}\overline{F}_{n}(x)=\overline{F}(x)
\end{equation*}

\noindent for all $x$ where $F(x)$ is continuous. From the expression of the density
of $F_{n}(x)$ we deduce that
\begin{equation*}
\overline{F}_{n}(x)=-\sum_{i=1}^{n}\frac{h_{n}(a_{i})}{a_{i}g_{n}'(a_{i})}\,e^{-a_{i}x},\
\qquad x>0\ .
\end{equation*}

\noindent Our objective is to prove that
\begin{equation}\label{DCT_main}
\lim_{n\to\infty}\sum_{i=1}^{n}\frac{h_{n}(a_{i})}{a_{i}g_{n}'(a_{i})}\,e^{-a_{i}x}=
\sum_{i=1}^{\infty}\frac{h(a_{i})}{a_{i}g'(a_{i})}\,e^{-a_{i}x},\
\qquad x>0\ .
\end{equation}

 The
proof of (\ref{DCT_main}) is done using the dominated convergence theorem. First notice that
for all $i\ge 1$,
\begin{equation}
\label{limitd}
\lim_{n\to\infty} g'_n(a_i)=g'(a_i).
\end{equation}
This is proved in the following way: On the one hand,
\begin{equation}
\label{limita}
\lim_{n\to\infty}g'_n(a_i)=-\frac{1}{a_i}\lim_{n\to\infty}
\prod_{\begin{smallmatrix}j=1\\j\ne i\end{smallmatrix}}^n
\Big(1-\frac{a_i}{a_j}\Big)=
-\frac{1}{a_i}
\prod_{\begin{smallmatrix}j=1\\j\ne i\end{smallmatrix}}^\infty
\Big(1-\frac{a_i}{a_j}\Big)
.
\end{equation}
On the other hand, for  $z\ne a_i,$
$$
\frac{g(z)}{1-z/a_i}=\prod_{\begin{smallmatrix}j=1\\j\ne i\end{smallmatrix}}^\infty
\Big(1-\frac{z}{a_j}\Big).
$$
The  function on the right  hand side is entire, so by analytic continuation,
\begin{equation}
\label{limitb}
\lim_{z\to a_i}\frac{g(z)}{1-z/a_i}=\prod_{\begin{smallmatrix}j=1\\j\ne i\end{smallmatrix}}^\infty
\Big(1-\frac{a_i}{a_j}\Big).
\end{equation}
Further, the function $1/g(z)$ has a simple pole at $a_i$. Hence,
\begin{equation}
\label{limitc}\lim_{z\to a_i}(z-a_i)\frac{1}{g(z)}=\text{Res}(1/g,a_i)=\frac{1}{g'(a_i)}.
\end{equation}
Combining (\ref{limita})--(\ref{limitc})   we obtain
(\ref{limitd}).

Fix $x>0$, set $\Lambda=\{a_{i},\, i\ge1\}$ and denote the counting measure by $\Gamma$.
By (\ref{limitd}),
\begin{equation*}
\lim_{n\to\infty}\frac{h_{n}(y)}{yg_{n}'(y)}\,e^{-yx}=\frac{h(y)}{yg'(y)}\,e^{-yx},
\qquad\forall\,y\in\Lambda\ .
\end{equation*}

\noindent Now observe that for $n>j$
$$\frac{h_{n+1}(a_j)}{g_{n+1}'(a_j)}=\frac{h_{n}(a_j)}{g_{n}'(a_j)}\
\frac{1-a_j/b_{n+1}}{1-a_j/a_{n+1}}\ ,$$
and since $a_{n+1}<b_{n+1}$ we conclude that
\begin{equation}\label{fita}
\left|\frac{h_{n}(a_j)}{g_{n}'(a_j)}\right|
\le \left|\frac{h_{n+1}(a_j)}{g_{n+1}'(a_j)}\right|
\le \left|\frac{h(a_j)}{g'(a_j)}\right|\ .
\end{equation}

\noindent Finally, let $\Phi(y)$ be the function
\begin{equation*}
\Phi(y):=\left|\frac{h(y)}{y{g}'(y)}\right|\1_{[0,1)}(y)e^{-yx}+
\left|\frac{h(y)}{{g}'(y)}\right|\1_{[1,\infty)}(y)e^{-yx},\
\qquad y>0\ .
\end{equation*}

\noindent Then
\begin{equation*}
\left|\frac{h_{n}(y)}{yg_{n}'(y)}e^{-yx}\right|
\le
\left|\frac{h_{n}(y)}{yg_{n}'(y)}\right|\1_{[0,1)}(y)e^{-yx}
+
\left|\frac{h_{n}(y)}{g_{n}'(y)}\right|\1_{[1,\infty)}(y)e^{-yx}\\
\le
\Phi(y),\ \qquad y>0\ .
\end{equation*}

\noindent By Lemma \ref{uniform_absolut},
\begin{equation*}
\int_{\Lambda}|\Phi|d\Gamma=
\sum_{i:\,a_{i}<1}\left|\frac{h(a_{i})}{a_{i}{g}'(a_{i})}\right|e^{-a_{i}x}
+\sum_{i:\,a_{i}\ge1}\left|\frac{h(a_{i})}{{g}'(a_{i})}\right|e^{-a_{i}x}
<\infty\ .
\end{equation*}

\noindent Hence we can apply the dominated convergence theorem to prove (\ref{DCT_main}). It
follows
\begin{equation*}
\overline{F}(x)=-\sum_{i=1}^{\infty}\frac{h(a_{i})}{a_{i}{g}'(a_{i})}\,e^{-a_{i}x},\ \qquad x>0\ .
\end{equation*}

\noindent Also from Lemma \ref{uniform_absolut}, we check that the set of continuity of
$F(x)$ is $(0,\infty)$. Moreover, from the uniform convergence on compact sets
of (\ref{general_D}),
\begin{equation*}
F'(x)=f(x)=-\sum_{n=1}^\infty \frac{h(a_n)}{g'(a_n)}\, e^{-a_n x},\quad x>0
\end{equation*}

\noindent as required. \qquad $\blacksquare$

\begin{remark}\label{main_h_1}
Standard manipulations show that Theorem \ref{main} is also true for $h\equiv 1$.
\end{remark}

\begin{remark}\label{density_laplace}
Notice that inequality (\ref{fita}) also holds true for Lemma \ref{tech_main} case (\ref{tech_2}).
Thus, Theorem \ref{main} for such $g(z)$ and $h(z)$ with orders lying in $(0,1)$ still holds
and concludes that the limiting density is
$$
f(x)=-\sum_{n=1}^\infty \frac{h(a_n)}{g'(a_n)}\, e^{-a_n |x|},\quad x\neq0\ .
$$

\noindent Moreover, the same proof of Theorem \ref{main} will be valid for functions $g(z)$
and $h(z)$ of order lying in $[1,2)$, provided that $\sigma_{a}\leq0$ for the series
(\ref{general_D}).
\end{remark}

The next result is an extension of the preceding theorem on the existence of a density
function. One would like to generalize such a result for characteristic functions
$\varphi(t)=h(it)/g(it)$ where $g(z)$ and $h(z)$ have zeros $\{\pm a_n,\, n\ge 1\}$ and
$\{\pm b_n,\, n\ge 1\}$ and orders lying in $[1,2)$. Entire functions of order greater than or equal
to one are much more difficult to treat than those that have  order less than one. Therefore we
will restrict the study to the setup of Lemma \ref{fundamental_lap}. As pointed out in the above remark,
we only need the absolute convergence of the general Dirichlet series to follow the same proof of
Theorem \ref{main} and obtain the density
\begin{equation}\label{series_lap}
f(x)=-\sum_{n=1}^\infty \frac{1}{g'(a_n)}\, e^{-a_n |x|},\quad x\neq0\ .
\end{equation}

\noindent However, we would like to present a different proof, Theorem \ref{main_lap}, which shows
that both problems, the exponential and Laplace convolutions, are the head and tail of the same
coin. It is worth  remarking that the random variables corresponding to that case are in the
homogeneous second Wiener chaos (see Janson \cite[chap. 6]{Jan97}), so our result is a step
 forward  in
the study of the densities of such an
 interesting space of random variables; in particular, this
case includes the L\'{e}vy area.

\begin{theorem}\label{main_lap}
Under assumptions of Lemma \ref{fundamental_lap}, the probability measure on $\R$
corresponding to the characteristic function $\varphi(t)=1/g(it)$ is absolutely continuous on
$\R\setminus\{0\}$ with (perhaps defective) density given by (\ref{series_lap}).
\end{theorem}

\noindent{\it Proof.}
Recall the expressions of $g(z)$ and $\wt g(z)$ defined in the proof of Lemma
\ref{fundamental_lap}. Notice that
\begin{equation}\label{aux_main_lap}
2a_{j}\wt g'(a_{j}^{2})=g'(a_{j})\ .
\end{equation}
Due to Proposition \ref{prop_existence} and Theorem \ref{main}, $\varphi(t)=1/\wt g(it)$ is a
characteristic function of a non--negative random variable, denoted by $Y$, with  density
$$f_Y(x)=\sum_{j=1}^\infty \frac{1}{\wt g'(a^2_j)}e^{-a^2_j  x},\quad x>0\ .$$
Fix $x\neq0$ and proceed heuristically as in Lemma \ref{lemma_lap_pos} using the afore mentioned
Bondesson argument and (\ref{aux_main_lap})  to obtain  the required expression
(\ref{series_lap}). To make the argument accurate we need to apply Fubini's theorem since
$f_{Y}(x)$ is an infinite series. It turns out that the absolute convergence of (\ref{general_D})
allows us to use Fubini's result. \qquad $\blacksquare$

\subsection{Infinite convolution of exponential and Laplace densities}

For the sake of completeness we will rewrite Remark \ref{main_h_1} and Theorem \ref{main_lap} in
a way that will extend Lemmas \ref{lemma_exp_pos} and \ref{lemma_lap_pos}. Moreover, we prove
that in this setup the resulting density is continuous in $[0,\infty)$ or $\R$.

Following Wintner, for a sequence of densities, $\{f_n,\, n\ge 1\}$, we will say that
$\star_{n=1}^\infty f_n$ is a convergent infinite convolution if the product
\begin{equation*}
\prod_{n=1}^{\infty}\psi_{n}(t)\ ,\qquad\text{ where }\qquad
\psi_{n}(t)=\int_{-\infty}^{\infty}e^{itx}f_{n}(x)dx\ ,
\end{equation*}

\noindent is uniformly convergent in every fixed finite $t$--interval.

\begin{proposition}\label{conv_exp}
Let $\{\lambda_{n}\, n\ge1\}$ be a strictly increasing sequence of positive numbers such that
for some $\rho\in(0,1)$, $\sum_{n\ge1}\lambda_{n}^{-\rho}<\infty$. Then
$\star_{n=1}^\infty {\mathcal E}{\rm xp}(\lambda_{n})$ converges to a  continuous density on
$[0,\infty)$
 which can be written as
\begin{equation}
\label{doble}
\star_{n=1}^\infty {\mathcal E}{\rm xp}(\lambda_{n})(x)=
\lim_{n\to\infty}\sum_{i=1}^{n}\frac{(-1)^{n+1}A(n)}{B(i,n)}\,e^{-\lambda_{i}x}=
\sum_{i=1}^{\infty}\lambda_i
\left[
\prod_{\begin{smallmatrix}k=1\\k\neq i\end{smallmatrix}}^{\infty}
    \left(1-\frac{\lambda_{i}}{\lambda_{k}}\right)\right]^{-1}
e^{-\lambda_{i}x}\ ,\ \ x>0,
\end{equation}
and
$$\star_{n=1}^\infty {\mathcal E}{\rm xp}(\lambda_{n})(0)=0.$$

\end{proposition}

\noindent{\it Proof.}
We will first point out that the infinite convolution is convergent. To this end we will use a
very useful result from Wintner \cite{Wintner35} that ensures the convergence of the infinite
convolution if
\begin{equation*}
\sum_{n=1}^{\infty}M_{n}<\infty\qquad\text{ where }\qquad
M_{n}=\E[|{\mathcal E}{\rm xp}(\lambda_{n})|]=
    \int_{-\infty}^{\infty}|x|\lambda_{n}\, e^{-\lambda_{n}x}\1_{[0,\infty)}dx\ .
\end{equation*}

\noindent Clearly $M_{n}=\lambda_{n}^{-1}$ and the condition to
guarantee the convergence of the infinite convolution is fulfilled. Moreover, if one
term of the infinite convolution has a continuous density  of bounded variation, then so has the
infinite convolution; and the continuous function $\star_{n=1}^m {\mathcal E}{\rm xp}(\lambda_{n})$
tends, as $m\to\infty$, to the infinite convolution uniformly in every bounded range (see Wintner
\cite{Wintner37}). It suffices to notice that
${\mathcal E}{\rm xp}(\lambda_{1})\star{\mathcal E}{\rm xp}(\lambda_{2})$ is continuous and of
 bounded
variation.
Hence, for all $x\ge 0$,
\begin{equation*}
\star_{n=1}^\infty {\mathcal E}{\rm xp}(\lambda_{n})(x)=
\lim_{n\to\infty}\sum_{i=1}^{n}\frac{(-1)^{n+1}A(n)}{B(i,n)}\, e^{-\lambda_{i}x}
\end{equation*}
To prove the second equality in (\ref{doble}), let
\begin{equation*}
g(z)=\prod_{k=1}^{\infty}\left(1- \frac{z}{\lambda_{k}}\right)\ .
\end{equation*}
\noindent Then $g(z)$ is an entire function of order less than $1$ and apply Theorem
\ref{main} to $x>0$, and the computations done in the first part of the proof of that theorem.
The value of the density  at zero is guarantied by Wintner  \cite{Wintner37} result
about the continuity on $\R$ of the density. \qquad $\blacksquare$

\bigskip

Note that the series in the right hand side of (\ref{doble}) may be divergent at $x=0$.


\begin{proposition}\label{conv_lap}
Let $\{\lambda_{n}\, n\ge1\}$ be a strictly increasing sequence of positive numbers such that
$\sum_{n\ge1}\lambda_{n}^{-2}<\infty$. Then $\star_{n=1}^\infty {\mathcal L}{\rm aplace}(\lambda_{n})$
converges to a continuous density on $\R$ which can be written as
\begin{equation*}
\star_{n=1}^\infty {\mathcal L}{\rm aplace}(\lambda_{n})(x)=
\lim_{n\to\infty}\sum_{i=1}^{n}\frac{(-1)^{n+1}A^{2}(n)}{E(i,n)}e^{-\lambda_{i}|x|}
\ ,\qquad x\in\R\ .
\end{equation*}
\end{proposition}

\noindent{\it Proof.}
As  in Proposition \ref{conv_exp}, we start by showing that the infinite convolution is
convergent. Jessen and Wintner show in \cite{JW35} that if
\begin{equation*}
\sum_{n=1}^{\infty}M_{n}^{1}\qquad\text{and}\qquad\sum_{n=1}^{\infty}M_{n}^{2}
\end{equation*}
are convergent, where $M_{n}^{1}$ and $M_{n}^{2}$ are the first and second moment of
${\mathcal L}{\rm aplace}(\lambda_{n})$ respectively, then so is the infinite convolution. One can
check that $M_{n}^{1}=0$ and $M_{n}^{2}=\lambda^{-2}_n$ to obtain the convergence of the infinite
convolution. Therefore we can apply Wintner \cite{Wintner37} to obtain the proposition since a
Laplace density is continuous and of finite variation. \qquad $\blacksquare$

\subsection{Existence of densities}\label{existence_d}

In this section we give three different conditions that guarantee  that the probability measure
corresponding to $h(it)/g(it)$ is absolutely continuous.

\begin{proposition}\label{no_atoms_1}
Under the hypotheses of Theorem \ref{main}, if $h(z)\equiv 1$, then the probability measure
corresponding to $1/g(it)$ is absolutely continuous and its density is continuous.
\end{proposition}

This result follows from Proposition \ref{conv_exp}.

\bigskip

In some cases it is known that the probability measure corresponding to $h(it)/g(it)$ has no atom at zero.
Since that probability measure is concentrated on $[0,\infty)$, this suffices for the existence of a density.

\begin{proposition}\label{no_atoms_2}
Assume the hypothesis of Theorem \ref{main} and denote by $\mu$ the probability measure
corresponding to $h(it)/g(it)$. If $\mu(\{0\})=0$, then $\mu$ is absolutely continuous.
\end{proposition}

The following lemma gives a sufficient  condition in order to apply a classical
criterion   for the existence of a continuous density.

\begin{lemma}
Let $g(z)$ and $h(z)$ be two entire functions  of order $\rho\in (0,1)$, both with non-zero
positive simple real zeros, and the zeros of $g(z)$ different from the zeros of $h(z)$. Denote
by $n_g(r)$ (respectively $n_h(r)$) the number  of the zeros of $g(z)$ with module less than $r$.
Assume the existence of the limits
$$\delta=\lim_{r\to\infty} \frac{n_g(r)}{r^\rho}>0 \qquad
\text{and} \qquad \delta'=\lim_{r\to\infty} \frac{n_h(r)}{r^\rho}>0,$$
with $\delta'<\delta$. Then
\begin{equation}\label{main_int}
\int_{-\infty}^{\infty}\left\vert \frac{h(it)}{g(it)}\right\vert\, dt <\infty\ .
\end{equation}

\end{lemma}

\noindent{\it Proof.}
Levin (\cite[p. 82, eqn. (2')]{Levin96}) proves that if $\theta\in (0,2\pi)$ and $r \to \infty$
then
\begin{equation}\label{asymp}
\log\vert g(re^{i\theta})\vert=
    \frac{\pi\delta r^\rho \cos\big(\rho(\theta-\pi)\big)}{\sin(\pi \rho)}
    +\frac{\ord(r^\rho)}{\sin(\theta/2)}\ .
\end{equation}

\noindent The analogous holds true for $h(re^{i\theta})$ with $\delta'$ instead of $\delta$.

Split integral (\ref{main_int}) into two parts as
\begin{equation}\label{suma-int}
\int_{-\infty}^{\infty}\Big\vert \frac{h(it)}{g(it)}\Big\vert\, dt
=\int_0^\infty \Big\vert \frac{h(te^{i\pi/2})}{g(te^{i\pi/2})}\Big\vert\, dt
+\int_0^\infty \Big\vert \frac{h(te^{i3\pi/2})}{g(te^{i3\pi/2})}\Big\vert\, dt.
\end{equation}

\noindent Due to (\ref{asymp}) there is $r_0>0$ and $C>0$ such that for $r>r_0$ we get
\begin{align*}
\left\vert \frac{h(te^{i\pi/2})}{g(te^{i\pi/2})}\right\vert
&=\exp\{-\cos(\rho\pi/2)\pi(\delta-\delta') \csc(\pi \rho) r^\rho+\ord(r^\rho)\}\\
&=\exp\{-r^\rho\big(\cos(\rho\pi/2)\pi(\delta-\delta')\csc(\pi\rho)+\ord(r^\rho)/r^\rho\big)\}
    \le e^{-Cr^\rho}.
\end{align*}

\noindent Hence
\begin{equation*}
\int_{r_0}^\infty \Big\vert \frac{h(te^{i\pi/2})}{g(te^{i\pi/2})}\Big\vert \, dt
\le
\int_{r_0}^\infty e^{-Cr^\rho}\, dr<\infty
\end{equation*}

\noindent since the last expression can be reduced to a convergent gamma integral. For the
integral over $[0,r_0]$ there is a straightforward bound using the same sort of derivations used
in Lemma \ref{fundamental}. The other integral on the right hand side of (\ref{suma-int}) is
bounded in a similar way. \qquad $\blacksquare$

With the notations of this lemma,

\begin{proposition}\label{no_atoms_3}
Under the hypothesis of Theorem \ref{main}, and assume the existence of the limits
$$\delta=\lim_{r\to\infty} \frac{n_g(r)}{r^\rho}>0 \qquad
    \text{and} \qquad \delta'=\lim_{r\to\infty} \frac{n_h(r)}{r^\rho}>0,$$
with $\delta'<\delta$. Then the probability measure corresponding to $h(it)/g(it)$ has a
continuous density.
\end{proposition}

\section{Examples}\label{examples}

We will now see different situations where we can apply the results obtained in the previous
sections. Most results are known, but here we get  them all using the same technique.

\subsection{L\'{e}vy area}

Let $\varphi(t)=\text{sech}(tT)$ be the characteristic function for the L\'{e}vy area (see
L\'{e}vy \cite{Levy51}), where the time component of the process varies in $[0,T]$. Here we set
$h\equiv1$ and $g(z)=\cos(zT)$, where $g(z)$ has order $1$. It is clear that
$\{\pm(2k-1)\pi/2T,\ k\ge1\}$ and $\{\pm(-1)^{k}/T,\ k\ge1\}$ are the poles and the residues of
$1/g(z)$ respectively, thus $\varphi(t)$ fulfills Proposition \ref{prop_existence_kent}.
Moreover, $g(z)$ is even and of exponential type, thus we apply Theorem \ref{main_lap} to obtain
\begin{equation*}
f(x)=\sum_{k=1}^{\infty}\frac{(-1)^{k+1}}{T}\,e^{-\frac{(2k-1)\pi}{2T}|x|}\qquad x\neq0\ .
\end{equation*}

\noindent Proposition \ref{conv_lap} ensures that the density function is continuous in $\R$.

\subsection{The first hitting time of a Bessel process}

The second example is a characteristic function obtained from
\begin{equation*}
k(z)=z^{-\nu/2}2^{\nu/2}\Gamma(\nu+1)J_{\nu}(\sqrt{2z})\ ,
\end{equation*}

\noindent where $J_{\nu}(z)$ is the Bessel function of first kind and order $\nu>-1$. Consider
$0<u<v<\infty$ and let $h(z)=k(u^2z)$ and $g(z)=k(v^2z)$. The probability measure of the corresponding
characteristic function $\varphi(t)=h(it)/g(it)$ describes the first hitting time of the point
$v$ by a Bessel process of order $\nu$ that starts at $u$, see Kent \cite{Kent78}, and can be
expressed as
\begin{equation*}
\varphi(t)=\left(\frac{v}{u}\right)^{\nu}\frac{J_{\nu}(u\sqrt{2it})}{J_{\nu}(v\sqrt{2it})}\ .
\end{equation*}

\noindent From the Taylor expansion of
\begin{equation*}
\left(\frac{2}{z}\right)^{\nu}J_{\nu}(z)
    =\sum_{n\ge0}\frac{(-1)^{n}}{n!\Gamma(n+\nu+1)4^n}z^{2n}\ ,
\end{equation*}
we can deduce the order of $k(z)$, since the order of an entire function is
\begin{equation*}
\rho=\limsup_{n\to\infty}\frac{n\ln(n)}{\ln(1/|c_{n}|)}\ ,
\end{equation*}
where $c_{n}$ are the coefficients of the Taylor expansion, see Levin \cite[p. 6]{Levin96}.
Due to the Stirling formula the above limit for $h(z)$ and $g(z)$ is $1/2$. Denote by
$\{j_{\nu,k},\, k\ge1\}$the positive zeros in order of magnitude of the Bessel function
$J_{\nu}(z)$ and by $a_{k}=j_{\nu,k}^2/(2v^2)$ the zeros of $g(z)$. Finally we can apply Theorem
\ref{main} to obtain the density
\begin{equation}\label{kent_dens}
f(x)=\sum_{k=1}^{\infty}
    \frac{j_{v,k}v^{\nu-2}J_{v}(j_{\nu,k}u/v)}{u^{\nu}J_{\nu+1}(j_{\nu,k})}
    e^{-\frac{j_{\nu,k}^{2}}{2v^2}x}\qquad x>0\ ,
\end{equation}

\noindent where we have used
\begin{equation*}
\frac{d}{dz}J_{\nu}(z)=-J_{\nu+1}(z)+\frac{\nu}{z}J_{\nu}(z)\ .
\end{equation*}

Equation (\ref{kent_dens}) is also derived by Borodin and Salminen \cite[p. 387]{BS02}. Notice
that the distribution function is absolutely continuous in $[0,\infty)$ since the probability
distribution gives no mass to $\{0\}$ due to the continuous paths of the Bessel process.

\subsubsection{Exit time from a $n$--dimensional sphere by a Brownian motion}

Let $T_{n}$ denote the random variable of the total time spent by an $n$--dimensional Brownian
motion starting at $0$ inside the sphere $S^{n-1}(r)$ of radius $r>0$ and $n\ge3$. Let $P_{n}$
denote the first exit time for a $n$--dimensional Brownian motion starting at $0$ from the sphere
$S^{n-1}(r)$ for $n\ge1$. Ciesielski and Taylor  \cite{CT62} show the remarkable equality of
the distribution functions of $T_{n}$ and $P_{n-2}$ for $n\ge3$, and they derive the
distribution function of $T_{n}$ using methods developed by Kac  \cite{Kac51}. They first
compute the solution for $n=3$ and then make a guess for the general framework. Finally they
compare the result with the distribution of $P_{n}$ which was computed by L\'{e}vy.

We can use Theorem \ref{main} to derive the density function of $T_{n}$ since Ciesielski and
Taylor give its characteristic function. In fact they establish the following result
\begin{equation*}
\E[e^{zT_{n}}]=\frac{(r\sqrt{2z})^{\nu-1}}{2^{\nu-1}\Gamma(\nu)J_{\nu-1}(r\sqrt{2z})}=
\prod_{i=1}^{\infty}\left(1-\frac{2r^{2}z}{j_{\nu-1,i}^{2}}\right)^{-1}\ ,
\end{equation*}

\noindent where $\nu=(n-2)/2$. This is a particular case of the previous example,
where we let $h(z)=1$ and $g(z)=k(r^2z)$. Use the same arguments as before to derive the density
function
\begin{equation*}
f(x)=\frac{1}{2^{\nu-1}\Gamma(\nu)r^{2}}\sum_{k=1}^{\infty}
\frac{j_{\nu-1,k}^{\nu}}{J_{\nu}(j_{\nu-1,k})}
e^{-\frac{j_{\nu-1,k}^{2}}{2r^{2}}x}\qquad x>0\ .
\end{equation*}

\noindent Notice that $T_{n}$ is absolutely continuous due to Proposition \ref{no_atoms_1}, and
we can set $f(0)=0$ to obtain a continuous density.

\subsubsection{The area under a squared Bessel bridge}

The characteristic function of $T_{n}$ is  easily representated when $\nu=n+1/2$ for $n\in\N$.
Let $\nu=3/2$ and $r=1$ to obtain the expression

\begin{equation*}
\varphi(t)=\frac{\sqrt{-2it}}{\sinh(\sqrt{-2it})}\ .
\end{equation*}

\noindent Revuz and Yor \cite[p. 465]{RevYor99} give the Laplace transform of the area
under a squared Bessel process starting at any point and arriving at zero. It turns out that the
above characteristic function corresponds to a Bessel process of order 2 starting and arriving at
 zero. Such functions appear recursively in the literature and many authors have studied them,
for instance \cite{BPY01} and \cite{PY03}.

The factorization of the characteristic function is particularly easy, and one can check the
identity
\begin{equation*}
\frac{1}{g(z)}=\prod_{k=1}^{\infty}\left(1-\frac{2z}{\pi^2 k^{2}}\right)^{-1}\ ,
\end{equation*}

\noindent where the residues of $1/g(z)$ are $\{(-1)^{k}\pi^2k^{2},\ k\geq 1\}$. Finally, the
density function can be written as
\begin{equation*}
f(x)=\left\{
\begin{array}{cr}
\sum_{k=1}^{\infty}(-1)^{k+1}\pi^2k^{2}e^{-\pi^2k^{2}x/2}& \text{for}\ x>0\\
0& \text{for}\ x=0
\end{array}
\right.\ ,
\end{equation*}
\noindent as stated in \cite{BPY01}. This characteristic function
also corresponds (with a change of parameters) to the square of a Kolmogorov law. The
corresponding distribution function is
\begin{equation*}
\vartheta_{4}=1+2\sum_{k=1}^{\infty}(-1)^{k}e^{-\pi^2k^{2}x/2}\qquad \text{for}\ x>0\ ,
\end{equation*}

\noindent which was obtained by Dugu\'{e} \cite{Dugue66}.

\subsection{Inverse Laplace transform}

Theta functions and related expressions have proved useful for manipulations of functionals of
Brownian motion. For instance, Borodin and Salminen use the inverse Laplace transform of
\begin{equation*}
\varphi(z)=\frac{v\sinh(u\sqrt{2z})}{u\sinh(v\sqrt{2z})}\quad 0<u<v\ ,
\end{equation*}

\noindent which turns out to be
\begin{equation}\label{example_laplace_1}
\mathcal{L}^{-1}(\varphi)(y)
=\frac{v}{u}\sum_{k=-\infty}^{\infty}\frac{v-u+2kv}{\sqrt{2\pi}y^{3/2}}
   \, e^{-\frac{(v-u+2kv)^{2}}{2y}}\quad y>0\ .
\end{equation}

\noindent Let us consider
$\varphi(t)=h(it)/g(it)$, where $g(z)=\frac{\sinh(v\sqrt{2z})}{v\sqrt{2z}}$ and similar for
$h(z)$. Standard manipulations lead to the computation of the sequences
$\{-\frac{k^{2}\pi^{2}}{2v^{2}},\ k\ge1\}$ and
$\{(-1)^{k+1}\frac{k\pi}{u}\sin\left(\frac{u}{v}k\pi\right),\ k\ge1\}$, which are the
poles and the residues of $h(z)/g(z)$ respectively. As pointed out in the previous example, both
functions $h(z)$ and $g(z)$ are entire functions of order $1/2$ and
$\varphi(t)$ satisfies the hypothesis of Proposition \ref{prop_existence}. Notice that the poles are negative and
hence we have to apply a generalization of Theorem \ref{main} for the negative case, this is
straightforward and we obtain the expression
\begin{equation*}
f(y)=\sum_{k=1}^{\infty}(-1)^{k+1}\frac{k\pi}{u}\sin\left(\frac{u}{v}k\pi\right)
e^{\frac{k^{2}\pi^{2}}{2v^{2}}y}\quad y<0\ .
\end{equation*}

\noindent Since we consider the Laplace transform we need to change the sign of $y$ in the above
expression and consider it in the range $(0,\infty)$. After doing that, it turns out that the
above series is equal to (\ref{example_laplace_1}) due to Poisson summation formulae. Moreover,
we are able to say that the density function $f(y)$ is continuous in $[0,\infty)$ due to
Proposition \ref{no_atoms_3} since
$n_{g}(r)=\left[\frac{v\sqrt{2r}}{\pi}\right]$, where $[x]$ stands for the integer part of $x$
and similar for $n_{h}(r)$.

\subsection{Heston density function}

The authors proved in \cite{BFU09} that the density function of the Heston model is
$\mathcal{C}^{\infty}$ and can be expressed as an infinite convolution of Bessel type densities.
We give here another expression in a particular case. In fact the search for such an expression
for the general case was the starting point for the present paper, thus we believe it is worth
recording it, even though it is only a partial result.

The Heston model for the log-spot is driven by the following system of stochastic differential
equations
\begin{equation*}
\left.
\begin{aligned}
dX_t    &=-\frac{1}{2} V_t\, dt+\sqrt{V_t}\, dZ_t\\
dV_t    &=a(b-V_t)\,  dt+c\sqrt{V_t}\, dW(t)
\end{aligned}
\right\}
\end{equation*}

\noindent where $a,b$ and $c$ are real positive constants. The processes $W$ and $Z$ are two
standard correlated Brownian motions
such that $\langle Z,W\rangle_t=\rho t$ for some $\rho\in[-1,1]$. For the particular case of
interest we set $2ab=c^2$ and consider  the volatility process $V$ to start at $0$. Then
the complex moment generating function of log--spot is
\begin{equation*}
\E[e^{zX_{t}}]=
\frac{e^{z(x_{0}-\rho ct/2)}e^{a}}{\cosh(P(z))+
\frac{a-\rho cz}{P(z)}
\sinh(P(z))}=\frac{e^{z(x_{0}-\rho ct/2)}}{g(z)}\ ,
\end{equation*}

\noindent where $x_{0}$ is the initial point for the process $X$ and
$P(z)=\sqrt{(a-c\rho z)^{2}+c^{2}(z-z^{2})}$. The term $e^{z(x_{0}-\rho ct/2)}$ is merely a
decentralisation term, while the function $g(z)$ is an entire function of order $1/2$ with all
zeros being real (see \cite{BFU09}). The function $g(z)$ has negative zeros for $\rho=-1$, while for
$\rho=1$ and $a\ge c$ all zeros are positive. In any case we can apply Theorem \ref{main} to
obtain an expression of the density function as a series of exponential type densities.
Obviously, the zeros of $g(z)$ must be computed numerically.


\begin{thebibliography}{10}

\bibitem{BPY01}
P.~Biane, J.~Pitman, and M.~Yor.
\newblock Probability laws related to the {J}acobi theta and {R}iemann zeta
  functions, and {B}rownian excursions.
\newblock {\em Bull. Amer. Math. Soc. (N.S.)}, 38(4):435--465 (electronic),
  2001.

\bibitem{Bondes92}
L.~Bondesson.
\newblock {\em Generalized gamma convolutions and related classes of
  distributions and densities}, volume~76 of {\em Lecture Notes in Statistics}.
\newblock Springer-Verlag, New York, 1992.

\bibitem{BS02}
A.~N. Borodin and P.~Salminen.
\newblock {\em Handbook of {B}rownian motion---facts and formulae}.
\newblock Probability and its Applications. Birkh\"auser Verlag, Basel, second
  edition, 2002.

\bibitem{CT62}
Z.~Ciesielski and S.~J. Taylor.
\newblock First passage times and sojourn times for {B}rownian motion in space
  and the exact {H}ausdorff measure of the sample path.
\newblock {\em Trans. Amer. Math. Soc.}, 103:434--450, 1962.

\bibitem{BFU09}
S.~del Ba\~{n}o~Rollin, A.~Ferreiro-Castilla, and F.~Utzet.
\newblock On the density of log--spot in {H}eston volatility model.
\newblock {\em Stochastic Process. Appl.},  120:2037--2063, 2010.

\bibitem{Dugue66}
D.~Dugu\'{e}.
\newblock Sur les lois de {K}olmogoroff et de van {M}ises.
\newblock {\em C. R. Acad. Sci. Paris S\'er. A-B}, 262:A999--A1000, 1966.

\bibitem{Fel66}
W.~Feller.
\newblock {\em An introduction to probability theory and its applications.
  {V}ol. {II}}.
\newblock John Wiley \& Sons Inc., New York, 1966.

\bibitem{Jan97}
S.~Janson.
\newblock {\em Gaussian {H}ilbert spaces}, volume 129 of {\em Cambridge Tracts
  in Mathematics}.
\newblock Cambridge University Press, Cambridge, 1997.

\bibitem{JW35}
B.~Jessen and A.~Wintner.
\newblock Distribution functions and the {R}iemann zeta function.
\newblock {\em Trans. Amer. Math. Soc.}, 38(1):48--88, 1935.

\bibitem{Kac51}
M.~Kac.
\newblock On some connections between probability theory and differential and
  integral equations.
\newblock In {\em Proceedings of the {S}econd {B}erkeley {S}ymposium on
  {M}athematical {S}tatistics and {P}robability, 1950}, pages 189--215,
  Berkeley and Los Angeles, 1951. University of California Press.

\bibitem{KT09}
M.~Katori and H.~Tanemura.
\newblock Zeros of airy function and relaxation process.
\newblock {\em Journal of Statistical Physics}, 136(6):1177--1204, September
  2009.

\bibitem{Kent78}
J.~Kent.
\newblock Some probabilistic properties of {B}essel functions.
\newblock {\em Ann. Probab.}, 6(5):760--770, 1978.

\bibitem{Levin96}
B.~Y. Levin.
\newblock {\em Lectures on entire functions}, volume 150 of {\em Translations
  of Mathematical Monographs}.
\newblock American Mathematical Society, Providence, RI, 1996.
\newblock In collaboration with and with a preface by Yu. Lyubarskii, M. Sodin
  and V. Tkachenko, Translated from the Russian manuscript by Tkachenko.

\bibitem{Levy51}
P.~L\'{e}vy.
\newblock Wiener's random function, and other {L}aplacian random functions.
\newblock In {\em Proceedings of the {S}econd {B}erkeley {S}ymposium on
  {M}athematical {S}tatistics and {P}robability, 1950}, pages 171--187,
  Berkeley and Los Angeles, 1951. University of California Press.

\bibitem{Luk70}
E.~Lukacs.
\newblock {\em Characteristic functions}.
\newblock Hafner Publishing Co., New York, 1970.
\newblock Second edition.

\bibitem{Man69}
S.~Mandelbrojt.
\newblock {\em S\'eries de {D}irichlet. {P}rincipes et m\'ethodes}, volume~11
  of {\em Monographies Internationales de Math\'ematiques Modernes}.
\newblock Gauthier-Villars, Paris, 1969.

\bibitem{PY03}
J.~Pitman and M.~Yor.
\newblock Infinitely divisible laws associated with hyperbolic functions.
\newblock {\em Canad. J. Math.}, 55(2):292--330, 2003.

\bibitem{RevYor99}
D.~Revuz and M.~Yor.
\newblock {\em Continuous martingales and {B}rownian motion}, volume 293 of
  {\em Grundlehren der Mathematischen Wissenschaften [Fundamental Principles of
  Mathematical Sciences]}.
\newblock Springer-Verlag, Berlin, third edition, 1999.

\bibitem{SchSonVon10}
R.~Schilling, L.~Song and  Z.~Vondra{\v{c}}ek,
\newblock {\em Bernstein functions, Theory and applications},
de Gruyter Studies in Mathematics, 37,
\newblock{Walter de Gruyter \& Co.,
Berlin, 2010.}




\bibitem{Tit52}
E.~C. Titchmarsh.
\newblock {\em The Theory of Functions}.
\newblock Oxford University press, second reprinted edition, 2002.

\bibitem{Wintner35}
A.~Wintner.
\newblock On the {D}ifferentiation of {I}nfinite {C}onvolutions.
\newblock {\em Amer. J. Math.}, 57(2):363--366, 1935.

\bibitem{Wintner37}
A.~Wintner.
\newblock On the {D}ensities of {I}nfinite {C}onvolutions.
\newblock {\em Amer. J. Math.}, 59(2):376--378, 1937.

\end{thebibliography}
\end{document}